\documentclass[12pt]{article}
\usepackage{amsmath,fullpage,amsthm,amssymb,xcolor,graphicx}

\newtheorem{thm}{Theorem}[section]
\newtheorem{corollary}{Corollary}[section]

\newtheorem{definition}{Definition}[section]
\newtheorem{remark}{Remark}[section]
\newtheorem{example}{Example}[section]

\title{On walk-regular graphs and optimal duals of frames generated by graphs}
\author{Deepshikha \thanks{Department of Mathematics, Shyampur Siddheswari Mahavidyalaya, University of Calcutta, West Bengal 711312, India. Email: dpmmehra@gmail.com }\   \and Aniruddha Samanta \thanks{Theoretical Statistics and Mathematics Unit, Indian Statistical Institute, Kolkata-700108, India. Email: aniruddha.sam@gmail.com}
}
\date{\today}
\begin{document}
	\maketitle
	\baselineskip=0.25in
	\begin{abstract}
Erasures are a common problem that arises while signals or data are being transmitted. A profound challenge in frame theory is to find the optimal dual frames ($OD$-frames) to minimize the reconstruction error if erasures occur. In this paper, we study the optimal duals of frames generated by graphs. First, we characterize walk-regular graphs. Then, it is shown that the diagonal entries of the Moore-Penrose inverse of the Laplacian matrix (or adjacency matrix) of a walk-regular graph are equal. Besides, we prove that connected graphs generate full spark frames. Using these results, we establish that the canonical dual frames are the unique $OD$-frames of a frame generated by a walk-regular graph.  A sufficient condition under which the canonical dual frame is the unique $OD$-frame is known. Here, we establish that the condition is also necessary if the frame is generated by a connected graph. 
	\end{abstract}

	{\bf AMS Subject Classification (2020):} 05C40, 05C50, 42C15, 42C40, 46C05.
	
	\textbf{Keywords.} Frames generated by graphs, Full spark frames, Optimal dual frames, Walk-regular graphs, Adjacency Matrix, Laplacian Matrix.
	
	\section{Introduction}
	
	Frame, graph, and matrix theories are the three applicable, useful, and fundamental branches of mathematics. They have their own deep significance in their domain of studies. In \cite{D}, the authors studied relations among these three different branches of mathematics by introducing the concept of frames generated by graphs. In this paper, we study the connection between graphs and optimal duals of frames generated by graphs. However, in order to establish a connection between graphs and optimal dual frames, we also explore some interesting properties of walk-regular graphs and their Moore-Penrose inverses. Before discussing the work presented in this paper, we briefly introduce spectral graph theory, frame theory, and optimal dual frames in the next three paragraphs. \vspace{4pt}

	Graph theory is one of the classical branches of mathematics. Graphs are used to model many real-life problems. The properties of graphs in terms of their associated matrix properties and vice-versa are studied in the field of spectral graph theory. Several matrices are associated with a graph, like adjacency matrix, Laplacian matrix, incidence matrix, distance matrix, etc. Along with the deep theoretical aspect of spectral graph theory, it has many real-life applications, which include measuring the connectivity of communication networks, hierarchical clustering, ranking of hyperlinks in search engines, mixing time of a random walk on graphs, image segmentation, and so on, see \cite{SGT2, SGT3, SGT4,SGT1}. 
	
	\vspace{4pt}
	
	In order to solve certain complex problems in nonharmonic Fourier series, Duffin and Schaeffer \cite{DS} developed the idea of frames. Daubechies, Grossmann, and Meyer \cite{DGM} demonstrated the beautiful significance of frame theory in data processing. A frame is a redundant sequence of vectors that provides infinitely many representations of elements of a vector space. Due to redundancy, frames have numerous applications like noise reduction, signal and image processing, sparse representations, data compression, etc., see \cite{BL, BP, CK, SH}.  In \cite{D}, one of the authors of this article introduced a novel class of finite frames generated by graphs. In this article, we draw connections between frame theory and spectral graph theory by studying the properties of frames generated by graphs.\vspace{4pt} 

		\vspace{4pt}
		
Erasures are a common issue in wireless sensor networks that arise from either transmission losses or disturbances during transmission. Overcoming the issue of erasures is crucial for the better recovery of signals. In order to address this issue, Lopez and Han proposed the idea of $OD$-frames in \cite{LH1}. The dual frames of a particular frame that can minimize the error caused by erasures are known as $OD$-frames, or optimal dual frames.  Many researchers introduced several measurements of error operators to investigate $OD$-frames,  see \cite{NAS, DS1, LH, LH1, PHM}. This article studies optimal duals of frames generated by graphs where the operator norm is used for the measurement of the error operator. A frame $\{f_i\}_{i\in[n]}$ for $\mathbb{C}^k$ is said to be a full spark frame if every subset of $\{f_i\}_{i\in[n]}$ with $k$ elements is linearly independent. Full spark frames do not appear very often, but they have a very special property that the loss or erasure of any $n-k$ frame elements still leaves a frame. In this article, we show that the frames generated by graphs are full spark frames.
	
		\vspace{4pt}
		
	The outline of this paper is as follows. In section 2, we provide the basics of frame theory, frames generated by graphs, optimal dual frames, and spectral graph theory required in the rest of the article. Section 3 studies walk-regular graphs. We mainly establish a characterization of a walk-regular graph. Then, we prove that the diagonal entries of the Moore-Penrose inverse of the Laplacian matrix (or adjacency matrix) of a walk-regular graph are equal. In section 4, we discuss about full spark frames. We prove that any frame generated by connected graphs is a full spark frame. Then, we find the spark of a frame generated by any graph. Results of Section 3 and Section 4 are used to establish the results of Section 5, where we study optimal duals of frames generated by graphs. We prove that the canonical dual frames of the frames generated by walk-regular graphs are unique $OD$-frames for any erasure. Further, a sufficient condition is shown for the canonical dual frame of a frame generated by a graph to be the $OD$-frame for 1-erasure. Theorem 2.6 in \cite{LH1} gives a sufficient condition under which the canonical duals are the unique $OD$-frames. In this article, we prove that the condition is also necessary if the frame is generated by a graph. We end the article by establishing that the canonical duals of frames generated by connected graphs are either the unique $OD$-frames or not $OD$-frames.	
	
	\section{Preliminaries}
In the rest of the article, every frame is a finite frame, $[n]:=\{1,2,\ldots,n\}$, and  $\mathbb{C}^k$  is a $k$-dimensional Hilbert space. We denote $M=[m_{i,j}]_{i\in [p], j\in [q]}$ as a matrix of order $p\times q$ where $m_{i,j}$ is the $(i,j)$-th entry of $M$. Also, $diag(\alpha_1,\alpha_2,\ldots,\alpha_n)$ denotes the diagonal matrix of order $n$ with respective diagonal entries $\alpha_1,\alpha_2,\ldots,\alpha_n$. 
	
	We begin with the following necessary definitions. 

	\begin{definition}
		A sequence $\{ f_i\}_{i\in[n]}$ of vectors in $\mathbb{C}^k$ is said to be a finite frame (or frame) for  $\mathbb{C}^k$ if there are positive constants $A_o\leq B_o$ such that for any $f\in\mathbb{C}^k$, we have
		\begin{align*}
			A_o\|f\|^2\leq\sum_{i\in[n]}|\langle f, f_i\rangle|^2\leq B_o\|f\|^2.
		\end{align*}
	\end{definition}
	
	The constants $B_o$ and $A_o$ are called an \emph{upper} and \emph{lower frame bound} of $\{f_i\}_{i\in[n]}$, respectively. 
	
\begin{definition}
			Frames $\{f_i\}_{i\in[n]}$ and $\{g_i\}_{i\in[n]}$ for $\mathbb{C}^k$ are called unitary equivalent frames if there is a unitary operator $U:\mathbb{C}^k\rightarrow\mathbb{C}^k$ such that $U(g_i)=f_i$ for $i\in[n]$.
	\end{definition}
	
	Suppose $F=\{f_i\}_{i\in[n]}$ is a frame for $\mathbb{C}^k$ then the \emph{analysis operator} $T:\mathbb{C}^k\rightarrow\mathbb{C}^n$ of $F$  is defined by  $$T(f)=\{\langle f,f_i\rangle\}_{i\in[n]}.$$
	The operator $T$ is a linear and bounded operator. The \emph{synthesis operator} of $F$ is the adjoint operator of $T$, defined as $$T^*(\{c_i\}_{i\in[n]})=\sum\limits_{i\in[n]}c_if_i.$$
	The canonical matrix representation of $T^*$ is a $k\times n$ matrix such that its $i^{\text{th}}$ column is the vector $f_i$ that is  $T^*=[f_1\,\,f_2\,\,\cdots\,\,f_n]$.

The composition of $T^*$ and $T$ is the \emph{frame operator} $S: \mathbb{C}^k\rightarrow \mathbb{C}^k$ defined as $$S(f) = T^*T(f) =\sum_{i\in[n]}\langle f, f_i\rangle f_i.$$
Here, $S$ is linear, bounded, invertible, self-adjoint, and positive. The frame operator provides the reconstruction formula as follows $$f = SS^{-1}f \ =\sum_{i\in[n]}\langle f, S^{-1}f_i\rangle f_i=\sum_{i\in[n]}\langle f, f_i \rangle S^{-1}f_i,  \ f\in\mathbb{C}^k.$$
	The sequence $\{S^{-1}f_i\}_{i\in[n]}$ is also a frame for $\mathbb{C}^k$ and it is known as the \emph{canonical dual frame} of $\{f_i\}_{i\in[n]}$. The values $\{\langle f,f_i\rangle\}_{i\in[n]}$ are known as the \emph{frame coefficients} of $f$ associated with the frame $\{f_i\}_{i\in[n]}$. The frame $\{h_i\}_{i\in[n]}$ is said to be a \emph{dual frame} or an \emph{alternate dual frame} of $\{f_i\}_{i\in[n]}$ if $$f =\sum_{i\in[n]}\langle f, f_i \rangle h_i \ =\sum_{i\in[n]}\langle f, h_i\rangle f_i, ~~\text{for any } f\in\mathbb{C}^k.$$
	
Another operator with respect to a frame  $\{f_i\}_{i\in[n]}$ is $TT^*:\mathbb{C}^n\rightarrow\mathbb{C}^n$, known as the \emph{Gramian operator}. Also, its canonical matrix representation is said to be the \emph{Gramian matrix}, defined by $\mathcal{G}=[\langle f_j,f_i\rangle]_{i,j\in[n]}$ where $\langle f_j,f_i\rangle$ is the $(i,j)$-th entry of $\mathcal{G}$. For more information about frame theory and its applications, we refer to \cite{BP, C, CK1, C2, CDH, HP}.

	\subsection{Optimal dual frames}
	If $F=\{f_i\}_{i\in[n]}$ is a frame for $\mathbb{C}^k$ and $H=\{h_i\}_{i\in[n]}$ is a dual frame of $F$, then for $\Lambda\subset[n]$, the error operator $E_{\Lambda,H}:\mathbb{C}^k\rightarrow\mathbb{C}^k$ is defined by 
	\begin{align*}
		E_{\Lambda,H}(f)=\sum_{i\in\Lambda}\langle f,f_i\rangle h_i.
	\end{align*}
	
	For $1\leq r<n$, let $D^r_{F,H}:=\max\{\|E_{\Lambda,H}\|:|\Lambda|=r\}$, where $\|E_{\Lambda,H}\|$ is the operator norm of $E_{\Lambda,H}$.
	So $D^r_{F,H}$ is the measurement of the maximum error if dual frame $H$ is used and $r$-erasures (i.e., loss of $r$ frame coefficients) occur. If $F$ is a frame with dual frame $\widetilde{H}$, then $\widetilde{H}$ is said to be an \emph{optimal dual frame} (or simply \emph{OD}-frame) for $1$-erasure if $$D^1_{F}:=\min\{D^1_{F,H}:H\text{ is a dual frame of }F\}=D^1_{F,\widetilde{H}}.$$
	 For $1<r<n$, if $\widetilde{H}$ is an $OD$-frame of $F$ for $(r-1)$-erasures and  $$D^r_{F}:=\min\{D^r_{F,H}:H\text{ is a dual frame of }F\}=D^r_{F,\widetilde{H}} $$ then we say that $\widetilde{H}$ is an $OD$-frame of $F$ for $r$-erasures. 
	 
	 Note that if $|\Lambda|=1$, then $D^1_{F,H}=\max\{\|E_{\Lambda,H}\|:|\Lambda|=1\}=\max\{\|f_i\|\|h_i\|:i\in[n]\}$. Let $F=\{f_i\}_{i\in[n]}$ be a frame. If $F$ has the frame operator $S$ and $c=\max\{\|f_i\|\|S^{-1}f_i\|:i\in[n]\}$, then we use the notation $\Lambda_{1, F}:=\{l:\|f_l\|\|S^{-1}f_l\|=c\}$.

In the following theorem, Lopez and Han \cite{LH1} presented a condition under which the canonical dual is the unique $OD$-frame. 
	\begin{thm}[\cite{LH1}]\label{thm2.1}
		Let $F=\{f_i\}_{i\in[n]}$ be a frame for a Hilbert space $\mathcal{H}$. If $F$ has the frame operator $S$ and $\|S^{-1}f_i\|\|f_i\|$ is constant for all $i\in[n]$, then the canonical dual is the unique $OD$-frame of $\{f_i\}_{i\in [n]}$ for any erasure.
	\end{thm}

\begin{thm}[\cite{LH}]\label{thm2.2}
		Let $F=\{f_i\}_{i\in[n]}$ be a frame for a Hilbert space $\mathcal{H}$. Suppose $S$ is the frame operator of $F$ and the set $\{f_i\}_{i\in\Lambda_{1, F}}$ is linearly independent. If there is a sequence of scalars $\{\alpha_i\}_{i\in[n]}$ such that $\sum\limits_{i\in[n]}\alpha_i f_i=0$ and $\alpha_i\neq 0$ for all $i\in\Lambda_{1, F}$, then $\{S^{-1}f_i\}_{i\in[n]}$ is not an $OD$-frame for $1$-erasure.
\end{thm}

	For the detailed study of optimal dual frames, readers may refer to \cite{DS1,LH,LH1,PHM}.

	\subsection{Matrices associated with graphs}
	Let $ G=(V(G),E(G)) $ be a simple graph with the vertex set $ V(G)=\{ v_1, v_2, \dots, v_n\}$  and edge set $ E(G) $. The vertices $ v_i $ and $ v_j $ are said to be adjacent if they are connected by an edge and is denoted by  $ v_i \sim v_j $. The number of vertices adjacent to $ v_i $ is known as the \emph{degree of the vertex} $ v_i $ and is denoted by $ d(v_i)$ or simply by $ d_i $. A \emph{walk} in a graph is a sequence of vertices and edges of the graph where both the vertices and edges in the sequence can be repeated. The number of edges in the sequence of a walk is called the \emph{length of the walk}. If the starting and ending vertices of a walk are identical, then it is known as a \emph{closed walk}. We denote a closed walk of length $p$ as a \emph{closed $p$-walk}. A graph $G$ is said to be a \emph{regular graph} if $d(v_i)=d(v_j)$ for all $i,j\in[n]$. If $d(v_i)=d(v_j)=r$ for all $i,j\in[n]$, then $G$ is called an \emph{$r$-regular graph}.  The transpose of a matrix $ M $ is denoted by $ M^t $. If $ M_1 $ and $ M_2 $ are matrices, then $ M_1 \oplus M_2 $ denotes the block matrix $\left[\begin{array}{cc}
		M_1 & \boldsymbol{0}  \\
		\boldsymbol{0} & M_2
	\end{array}\right]$.

The degree matrix of a graph $G$ on $n$ vertices, denoted as $D(G):=diag(d_1,d_2,\ldots,d_n)$. The \emph{adjacency matrix} of $ G $ on $ n $ vertices is $ A(G):=(a_{i,j})_{n\times n} $, where
	\begin{align*}
	a_{i,j}=\begin{cases}
		1, \ \text{for }  v_i\sim v_j \text{ and } i\neq j \\
		0, \text{ elsewhere.}
	\end{cases}
\end{align*}
The matrix $ A(G) $ is simply written as $ A $, if the graph $G$ is clear from the context. The \emph{eigenvalues of a graph $ G $} are the eigenvalues of $ A(G) $. The Laplacian matrix of a graph $G$ is defined as $L(G):=D(G)-A(G)$. We simply write $ L(G) $ as $ L $ if $ G $ is clear from the context.  The Laplacian matrix is a positive semi-definite matrix. For a simple graph $G$ on $n$ vertices with $p$ components, the rank of $L(G)$ is $n-p$. Readers may refer to \cite{Bapat, DMG, DMG1} for more information about adjacency and Laplacian matrices.


Let us recall the notion of walk-regular graphs, an important class of graphs in this article.

\begin{definition}[\cite{GM}]
	A simple graph $ G $ is called a walk-regular graph if, for any $ k\in \mathbb{N}$, the number of closed $k$-walks beginning at the vertex $ v $ is independent of the choice of the vertex $ v $.
\end{definition}

\begin{thm}[\cite{GM}]
	Suppose $G$ is a simple graph with adjacency matrix $A$. Then, $G$ is walk-regular if and only if for any $i\in\mathbb{N}$, the diagonal entries of $A^i$ are equal.
\end{thm}

Next, we see the definition of the Moore-Penrose inverse of a matrix.
\begin{definition}[\cite{Bapat}]
	If $B$ is a square matrix, then the Moore-Penrose inverse of $B$ is denoted by $B^+$, such that $BB^+B=B$, $B^+BB^+=B^+$, $(BB^+)^t=BB^+$, and $(B^+B)^t=B^+B$.
\end{definition}
Note that the Moore-Penrose inverse of every square matrix exists and it is unique.
\begin{thm}[\cite{Bapat}]\label{th2.4}
Let $B$ be a symmetric matrix of order $n$. If $M$ is an orthogonal matrix such that $B=M\,diag(\lambda_1,\lambda_2,\ldots,\lambda_n)\,M^*$ where $\lambda_1,\lambda_2,\ldots,\lambda_k$ are the non-zero eigenvalues of $B$, then $B^+=M\,diag\left(\frac{1}{\lambda_1},\frac{1}{\lambda_2},\ldots,\frac{1}{\lambda_k}, 0, \dots, 0\right)\, M^*$.
\end{thm}

	\subsection{Frames generated by graphs}
	
	In \cite{D}, one of the authors in the present paper introduced the concept of $L_G(n,k)$-frames and frames generated by graphs. First, let us see the definition of $L_G(n,k)$-frames.
	\begin{definition}
		Let $G$ be a simple graph on  $n$ vertices and has $n-k$ components. If $L$ is the Laplacian matrix of $G$ such that $L=M\, diag(\lambda_1,\lambda_2,\ldots,\lambda_{k},0,\ldots,0)\,M^*$. If $\{e_i\}_{i\in[n]}$ is the standard canonical orthonormal basis of $\mathbb{C}^n$ and $B=diag\left(\sqrt{\lambda_1},\ldots,\sqrt{\lambda_k}\right)M_1^*$, where $M_1$ is a submatrix of $M$ formed by the first $k$ columns, then $\{B(e_i)\}_{i\in[n]}$ is called an $L_G(n,k)$-frame for $\mathbb{C}^k$.
	\end{definition}
	Following is the definition of frames generated by graphs, which are also called $G(n,k)$-frames.
	\begin{definition}
		Let $\{f_i\}_{i\in [n]}$ be a frame for $\mathbb{C}^k$ with Gramian matrix $\mathcal{G}$. If $G$ is a graph with Laplacian matrix $L$ such that $\mathcal{G}=L$, then $\{f_i\}_{i\in [n]}$ is called a frame generated by the graph $G$ for $\mathbb{C}^k$. In short, we call $\{f_i\}_{i\in [n]}$ as a $G(n,k)$-frame for $\mathbb{C}^k$.
	\end{definition}
	In \cite{D}, it is shown that every $L_G(n,k)$-frame is a $G(n,k)$-frame. The following result shows that any two frames generated by the same graph are unitary equivalent.
	
	\begin{thm}[\cite{D}]\label{thm2.6}
		If $G$ is a graph, then frames generated by $G$ are unitary equivalent.
	\end{thm}
	In the next theorem, it is shown that the frame operator of any $L_G(n,k)$-frame is a diagonal matrix. 
	\begin{thm}[\cite{D}]\label{thm2.7}
		If $G$ is a graph and $F$ is an $L_G(n,k)$-frame, then the frame operator $S$ of the frame $F$ is a diagonal matrix. Further, if $\lambda_1,\lambda_2,\ldots,\lambda_k$ are the non-zero eigenvalues of the Laplacian matrix of $G$, then $S=diag(\lambda_1,\lambda_2,\ldots,\lambda_k)$.
	\end{thm}
	The following theorem gives the family of dual frames of the frames generated by graphs.
	\begin{thm}\label{thm2.8}
		Let $G$ be a simple graph with components $G_1,G_2,\ldots, G_{m}$ having vertex sets $\{v_1,v_2,\ldots,v_{n_1}\}, \{v_{n_1+1}, v_{n_1+2}, \ldots, v_{n_2}\}, \ldots,\{v_{n_{m-1}+1}, v_{n_{m-1}+2}, \ldots, v_{n_m}(=v_n)\}$, respectively. If $F=\{f_i\}_{i\in[n]}$ is a $G(n,n-m)$-frame with the frame operator $S$, then any dual frame of $F$ is of the form $\{S^{-1}f_i+\nu_1\}_{i=1}^{n_1}\bigcup\{S^{-1}f_i+\nu_2\}_{i=n_1+1}^{n_2}\bigcup\cdots\bigcup\{S^{-1}f_i+\nu_m\}_{i=n_{m-1}+1}^{n}$ where $\nu_1,\nu_2,\ldots,\nu_m$ are arbitrary vectors in $\mathbb{C}^{n-m}$.
	\end{thm}
%
	
	\section{On walk-regular graphs}
	Our objective in this section is to study walk-regular graphs and their Moore-Penrose inverses, which will be useful in the later sections.  We begin by computing the determinant of the following matrices in order to characterize the walk-regular graphs.
	
	\vspace{8pt}
	Suppose $M=\left[\begin{array}{cc}
		a_1 & a_2  \\[1.5mm]
		a_1^2 & a_2^2
	\end{array}\right]$. Then, $ det(M)=\left|\begin{array}{cc}
	a_1 & a_2  \\[1.5mm]
	a_1^2 & a_2^2
\end{array}\right|=a_1a_2(a_2-a_1). $

\vspace{8pt}
Now suppose $M=\left[\begin{array}{ccc}
	a_1 & a_2  & a_3\\[1.5mm]
	a_1^2 & a_2^2 & a_3^2\\[1.5mm]
	a_1^3 & a_2^3 & a_3^3
\end{array}\right]$. Then, by inductive steps, we have
\begin{align*}
	det(M)&=\left|\begin{array}{ccc}
		a_1 & a_2  & a_3\\[1.5mm]
		a_1^2 & a_2^2 & a_3^2\\[1.5mm]
		a_1^3 & a_2^3 & a_3^3
	\end{array}\right|=a_1(a_2-a_1)(a_3-a_1)\left|\begin{array}{cc}
	a_2 & a_3\\[1.5mm]
	a_2^2 & a_3^2
\end{array}\right|=a_1a_2a_3(a_2-a_1)(a_3-a_1)(a_3-a_2).
\end{align*}
Therefore, continuing the above process, we have the following remark.

\begin{remark}\label{rem3.1}
If
$$M=\left[\begin{array}{cccc}
	a_1 & a_2  & \cdots & a_n\\[1.5mm]
	a_1^2 & a_2^2 & \cdots & a_n^2\\[1.5mm]
	\vdots & \vdots & \ddots & \vdots\\[1.5mm]
	a_1^n & a_2^n & \cdots & a_n^n
\end{array}\right]$$
then $det(A)=a_1a_2\ldots a_n\prod\limits_{i,j\in[n],i>j}(a_i-a_j)$.
\end{remark}

A simple graph is a walk-regular graph if, for any $p\in \mathbb{N}$, the number of closed $p$-walks at any vertex $v$ is independent of the choice of the vertex $v$. In the following theorem, we show that it is not necessary to check closed walks of all lengths for walk-regular graphs.
\begin{thm}\label{thm3.1}
Let $ G $ be a simple graph with $ k $ distinct non-zero eigenvalues. For $ p\in [k] $, if the number of closed $p$-walk is independent of the choice of vertices, then $ G $ is walk-regular.	
\end{thm}
\proof
Let $ A$ be the adjacency matrix of $ G $ with non-zero eigenvalues $\lambda_1,\lambda_2,\ldots,\lambda_s$. Without loss of generality, assume $\lambda_1,\lambda_2,\ldots,\lambda_k$ are distinct non-zero eigenvalues of $A$. Let $A=MDM^*$ such that $M=[m_{i,j}]_{i,j\in[n]}$ and $D=diag(\lambda_1,\lambda_2,\ldots,\lambda_s,0,\ldots,0)$. For $p \in [k]$, the number of closed $p$-walks is independent of choice of vertices, so $ A^p $ has equal diagonal entries. Consider the first and second diagonal entries of $A^p$, as they are equal, thus $\lambda_1^p|m_{1,1}|^2+\cdots+\lambda_s^p|m_{1,s}|^2=\lambda_1^p|m_{2,1}|^2+\cdots+\lambda_s^p|m_{2,s}|^2$. After combining equal eigenvalues, we get $\lambda_1^p|\widetilde{m_{1,1}}|^2+\cdots+\lambda_k^p|\widetilde{m_{1,k}}|^2=\lambda_1^p|\widetilde{m_{2,1}}|^2+\cdots+\lambda_k^p|\widetilde{m_{2,k}}|^2$ that is $\lambda_1^p(|\widetilde{m_{1,1}}|^2-|\widetilde{m_{2,1}}|^2)+\cdots+\lambda_k^p(|\widetilde{m_{1,k}}|^2-|\widetilde{m_{2,k}}|^2)=0$, for all $p \in [k] $. Let $|\widetilde{m_{1,l}}|^2-|\widetilde{m_{2,l}}|^2=x_l$ for $l\in [k]$. Let $$B=\left[\begin{array}{cccc}
	\lambda_1 & \lambda_2  & \cdots & \lambda_k\\
	\lambda_1^2 & \lambda_2^2 & \cdots & \lambda_k^2\\
	\vdots & \vdots & \ddots & \vdots\\
	\lambda_1^k & \lambda_2^k & \cdots & \lambda_k^k
\end{array}\right], X=\left[\begin{array}{c}
	x_1\\
	x_2\\
	\vdots \\
	x_k
\end{array}\right].$$ Then, we have $BX=0$. By Remark \ref{rem3.1}, $det(B)=\lambda_1\lambda_2\cdots \lambda_k\prod\limits_{i,j\in[k],i>j}(\lambda_i-\lambda_j)\neq 0$ as $\lambda_1,\lambda_2,\ldots,\lambda_k$ are distinct. Thus, $X=0$. Hence, $x_l=0$ for all $l\in[k]$. Now, for any $i\in\mathbb{N}$, $\lambda_1^ix_1+\cdots+\lambda_k^ix_k=0$ that is $\lambda_1^i|\widetilde{m_{1,1}}|^2+\cdots+\lambda_k^i|\widetilde{m_{1,k}}|^2=\lambda_1^i|\widetilde{m_{2,1}}|^2+\cdots+\lambda_k^i|\widetilde{m_{2,k}}|^2$. Thus, the first and second diagonal entries of $ A^i $ are equal for $ i \in \mathbb{N} $. Proceeding in this way, we can show that the diagonal entries of $A^i$ are equal for $i\in\mathbb{N}$. Therefore, the number of closed $i$-walks is independent of choice of vertices for any $ i \in \mathbb{N} $, so $ G $ is walk-regular. 
\endproof

The following result provides a characterization of the walk-regular graphs, which can be obtained from the above theorem.
\begin{corollary}
		Suppose $G$ is a simple graph with $k$ distinct non-zero eigenvalues. Then, $G$ is walk-regular if and only if the number of closed $p$-walks is independent of the choiceof vertices for any $ p \in [k] $. 
\end{corollary}

Next, we study about the diagonal entries of the Moore-Penrose inverse of matrices associated with the walk-regular graphs.
\begin{thm}\label{thm3.2}
If $ G $ is a walk-regular graph with adjacency matrix $ A $, then the Moore-Penrose inverse $ A $ has equal diagonal entries.	
\end{thm}
\proof
Let $\lambda_1,\lambda_2,\ldots,\lambda_s$ be the non-zero eigenvalues of $A$ in which $\lambda_1,\lambda_2,\ldots,\lambda_k$ are distinct. If $A=MDM^*$ such that the matrices $M=[m_{i,j}]_{i,j\in[n]}$ and $D=diag(\lambda_1,\lambda_2,\ldots,\lambda_s,0,\ldots,0)$, then by Theorem \ref{th2.4}, the Moore-Penrose inverse $A^+=MD_1M^*$ where $D_1=diag(\frac{1}{\lambda_1},\frac{1}{\lambda_2},\ldots,\break \frac{1}{\lambda_s}, 0,\ldots,0)$. Let $X=\left[\begin{array}{c}
	x_1\\
	x_2\\
	\vdots\\
	x_k
\end{array}\right]$ be same as in the proof of Theorem \ref{thm3.1}. Then, $x_1=x_2=\cdots=x_k=0$. Thus $\frac{1}{\lambda_1}x_1+\frac{1}{\lambda_2}x_2+\cdots+\frac{1}{\lambda_k}x_k=0$. Hence, $\frac{1}{\lambda_1}|m_{1,1}|^2+\cdots+\frac{1}{\lambda_m}|m_{1,m}|^2=\frac{1}{\lambda_1}|m_{2,1}|^2+\cdots+\frac{1}{\lambda_m}|m_{2,m}|^2$. Therefore, the first and second diagonal entries of $ A^+ $ are equal. Similarly, one can show that all the diagonal entries of $A^+$ are equal.
\endproof

\begin{thm}\label{thm3.3}
If $L$ is the Laplacian matrix of a walk-regular graph $G$, then the Moore-Penrose inverse $L^+$ of $L$ has equal diagonal entries.
\end{thm}
\proof
Let $G$ be a walk-regular graph. Then, $G$ is regular, say $r$-regular. Let $A$ be the adjacency matrix of $G$. Therefore, the diagonal entries of $A^k$ are equal for $k\in\mathbb{N}$. Since $L=rI-A$, so for $k\in\mathbb{N}$, $L^k$ has equal diagonal entries. Now, by using the similar arguments as used in Theorem \ref{thm3.2}, we obtain that $L^+$ has equal diagonal entries. 
\endproof

	\section{Full Spark frame}\label{sec1}
	The section starts with the definition of full spark frames, which is an important class of frames.
	\begin{definition}[\cite{CK1}]
		Suppose $\{f_i\}_{i\in[n]}$ is a frame for $\mathbb{C}^k$ then the cardinality of the smallest linearly dependent subset of $\{f_i\}_{i\in[n]}$ is called the spark of the frame. If spark is $k+1$, then we say that $\{f_i\}_{i\in[n]}$ is a full spark frame for $\mathbb{C}^k$.
	\end{definition}
	Full spark frames have the special feature that the loss of some frame elements still leaves a frame. In general, it is difficult to find full spark frames. The following theorem shows that for a connected graph $G$, $L_G(n,k)$-frames are full spark frames. 
	\begin{thm}\label{thm4.1}
		If $G$ is a connected graph and $F=\{f_i\}_{i\in[n]}$ is an $L_G(n,n-1)$-frame, then $F$ is a full spark frame.
	\end{thm}
	\proof
	Suppose $\{ v_1,v_2, \dots, v_n\}$ is the vertex set of $G$ and $\{e_i\}_{i\in [n]}$ is the standard canonical orthonormal basis of $\mathbb{C}^n$. Let $L$ be the Laplacian matrix of $G$ with eigenvalues $\lambda_1\geq\lambda_2\geq\ldots\geq\lambda_{n-1}>\lambda_n=0$. Note that the rank of $L$ is $n-1$. Then, there exist an orthogonal matrix $M=[m_{i,j}]_{i,j \in [n]}$ consisting of eigenvectors of $L$ and diagonal matrix $D=diag(\lambda_1,\lambda_2,\ldots,\lambda_{n-1},0)$ such that $L=MDM^*$ and  $f_i=Be_i$ for $i\in [n]$, where $B=D_1M_1^*$, $M_1=[m_{i,j}]_{i\in [n], j\in [n-1]}$, and $D_1=diag(\sqrt{\lambda_1},\sqrt{\lambda_2},\ldots,\sqrt{\lambda_{n-1}})$. 
	
	 Let $m_j=[m_{i,j}]_{i \in [n]}$ for $j\in[n]$. Then, $m_j$ is an eigenvector of $L$ corresponding to the eigenvalue $\lambda_j$. Thus $Lm_j=\lambda_jm_j$ that is 
	\begin{align*}
		\left[
		\begin{array}{c}
			d_1m_{1,j}-\sum\limits_{i\sim 1}m_{i,j}\\
			\vdots\\
			d_nm_{n,j}-\sum\limits_{i\sim n}m_{i,j}
		\end{array}
		\right]
		=	\left[
		\begin{array}{c}
			\lambda_jm_{1,j}\\
			\vdots\\
			\lambda_jm_{n,j}
		\end{array}
		\right].
	\end{align*}
	Thus, we have
	\begin{align*}
		\lambda_j(m_{1,j}+\cdots+m_{n,j})&=d_1m_{1,j}-\sum\limits_{i\sim 1}m_{i,j}+\cdots+d_nm_{n,j}-\sum\limits_{i\sim n}m_{i,j}=0.
	\end{align*}
	Hence, $m_{1,j}+\cdots+m_{n,j}=0$. Then, we have
	\begin{align*}
		B=D_1M_1^*=\left[\begin{array}{c}
			\sqrt{\lambda_1}m_1^*\\
			\vdots\\
			\sqrt{\lambda_{n-1}}m_{n-1}^*
		\end{array}\right]
		=\left[\begin{array}{cccc}
			\sqrt{\lambda_1}\overline{m_{1,1}} & \sqrt{\lambda_1}\overline{m_{2,1}} & \cdots & \sqrt{\lambda_1}\overline{m_{n,1}}\\
			\vdots & \vdots & \ddots & \vdots\\
			\sqrt{\lambda_{n-1}}\overline{m_{1,n-1}} &  \sqrt{\lambda_{n-1}}\overline{m_{2,n-1}} & \cdots & \sqrt{\lambda_{n-1}}\overline{m_{n,n-1}}
		\end{array}\right].
	\end{align*}   
 \newline
Thus, each row sum of the matrix $B$ is $0$. That is $f_1+f_2+\cdots+f_n=0$. Thus, for any $i\in [n]$, $f_i=-\sum\limits_{j\neq i}f_j$. Hence, $f_i\in span\{f_j\}_{j\in[n],j\neq i}$ and this gives $\mathbb{C}^{n-1}=span\{f_1,\ldots,f_n\}\subseteq span\{f_j\}_{j\in [n],j\neq i}$. Thus, $\{f_j\}_{j\in [n],j\neq i}$ is linearly independent. Therefore,  $F$ is a full spark frame.
	\endproof

	The next theorem shows that any frame generated by a connected graph is a full spark frame.
	
	\begin{thm}\label{thm4.2}
		If $G$ is a connected graph of $n$ vertices and $F=\{f_i\}_{i\in[n]}$ is a $G(n,n-1)$-frame for $\mathbb{C}^{n-1}$, then $F$ is a full spark frame.
	\end{thm}
	\proof
	Let $\widetilde{F}=\{\widetilde{f}_i\}_{i\in[n]}$ be an $L_{G}(n,n-1)$-frame for $\mathbb{C}^{n-1}$. Then, by Theorem \ref{thm4.1}, $\widetilde{F}$ is a full spark frame. Also, by Theorem \ref{thm2.6}, $F$ and $\widetilde{F}$ are unitary equivalent frames. Thus, there exist a unitary operator $U:\mathbb{C}^{n-1}\rightarrow\mathbb{C}^{n-1}$ such that $U(\widetilde{f}_i)=f_i$ for $i\in[n]$. Now for any $\{m_1,m_2,\ldots,m_{n-1}\}\subset [n]$, $\{\widetilde{f}_{m_1},\widetilde{f}_{m_2},\ldots,\widetilde{f}_{m_{n-1}}\}$ is linearly independent. Hence, $\{U(\widetilde{f}_{m_1}),U(\widetilde{f}_{m_2}),\ldots,U(\widetilde{f}_{m_{n-1}})\}=\{f_{m_1},f_{m_2},\ldots,f_{m_{n-1}}\}$ is linearly independent. Therefore, $F$ is a full spark frame.
	\endproof
	
	Theorem \ref{thm4.1} and Theorem \ref{thm4.2} need not hold for frames generated by disconnected graphs is illustrated in the following example.
	\begin{example}\label{Exa4.1}
	{\em 	Consider the disconnected graph $G$  given in the Figure \ref{fig1}. Then
		\begin{figure}
			\begin{center}
				\includegraphics[scale= 0.65]{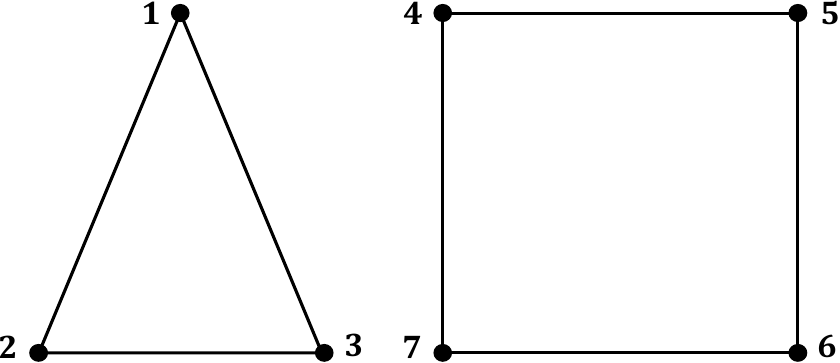}
				\caption{ Graph $ G $} \label{fig1}
			\end{center}
		\end{figure}
  
		\begin{align*}
			L(G)=\left[\begin{array}{ccccccc}
				2 & -1 & -1 & 0 & 0 & 0 & 0\\
				-1 & 2 & -1 & 0 & 0 & 0 & 0\\
				-1 & -1 & 2 & 0 & 0 & 0 & 0\\
				0 & 0 & 0 & 2 &-1 & 0 &-1\\
				0 & 0 & 0 & -1 & 2 & -1 & 0\\
				0 & 0 & 0 & 0 &-1 & 2 &-1\\
				0 & 0 & 0 & -1 & 0 & -1 & 2
			\end{array}
			\right]=MDM^*
		\end{align*}
		where $D=diag(4,2,2,3,3,0,0)$ and $M=\left[\begin{array}{ccccccc}
			0 & 0 & 0 & \frac{1}{\sqrt{2}} & \frac{-1}{\sqrt{6}} & \frac{1}{\sqrt{3}} & 0\\
			0 & 0 & 0 & 0 & \frac{2}{\sqrt{6}} & \frac{1}{\sqrt{3}} & 0\\
			0 & 0 & 0 & \frac{-1}{\sqrt{2}} & \frac{-1}{\sqrt{6}} & \frac{1}{\sqrt{3}} & 0\\
			\frac{1}{2} & \frac{1}{\sqrt{2}} & 0 & 0 & 0 & 0 & \frac{1}{2}\\
			\frac{-1}{2} & 0 & \frac{1}{\sqrt{2}} & 0 & 0 & 0 & \frac{1}{2}\\
			\frac{1}{2} & \frac{-1}{\sqrt{2}} & 0 & 0 & 0 & 0 & \frac{1}{2}\\
			\frac{-1}{2} & 0 & \frac{-1}{\sqrt{2}} & 0 & 0 & 0 & \frac{1}{2}
		\end{array}
		\right]$. Let $F$ be an $L_G(7,5)$-frame for $\mathbb{C}^5$ such that the synthesis operator of $F$ is $diag(2,\sqrt{2},\sqrt{2},\sqrt{3},\sqrt{3})M_1^*$ where $M_1$ is the matrix obtained from the matrix $M$ by taking the first five columns of $M$. Then, we have
		\begin{align*}
			F=\{f_i\}_{i\in [7]}=\left\{\left[\begin{array}{c}
				0\\
				0\\
				0\\
				\frac{\sqrt{3}}{\sqrt{2}}\\
				\frac{-\sqrt{3}}{\sqrt{6}}
			\end{array}\right],\left[\begin{array}{c}
				0\\
				0\\
				0\\
				0\\
				\frac{2\sqrt{3}
				}{\sqrt{6}}
			\end{array}\right], \left[\begin{array}{c}
				0\\
				0\\
				0\\
				\frac{-\sqrt{3}}{\sqrt{2}}\\
				\frac{-\sqrt{3}}{\sqrt{6}}
			\end{array}\right], \left[\begin{array}{c}
				1\\
				1\\
				0\\
				0\\
				0
			\end{array}\right], \left[\begin{array}{c}
				-1\\
				0\\
				1\\
				0\\
				0
			\end{array}\right], \left[\begin{array}{c}
				1\\
				-1\\
				0\\
				0\\
				0
			\end{array}\right], \left[\begin{array}{c}
				-1\\
				0\\
				-1\\
				0\\
				0
			\end{array}\right]\right\}
		\end{align*}
		Now $\{f_1,f_2,f_3,f_4,f_5\}$ is a subset of $F$ with five vectors but it is linearly dependent. Hence, $F$ is not a full spark frame.}
	\end{example}
	
	In Theorem \ref{thm4.2}, we see that the frames generated by connected graphs are full spark frames. However, the frames generated by disconnected graphs need not be full spark, which is justified in Example \ref{Exa4.1}. In the following theorem, we establish a relation between the spark of a frame and the vertex cardinality of a graph.
	\begin{thm}\label{thm4.3}
		Let $G$ be a graph of $n$ vertices with $k$ components and $F=\{f_i\}_{i\in[n]}$ be an $L_G(n,n-k)$-frame for $\mathbb{C}^{n-k}$. If $ m $ is the minimum cardinality of the vertex set of components of $ G $, then the spark of $F$ is $m$.
	\end{thm}
	\proof
Suppose $G_1,\ldots,G_{k}$ are the components of $G$. Take $ n_i=|V(G_i)| $, for $ i\in [k] $.  Then, $L(G)=L(G_1)\oplus\cdots\oplus L(G_{k})$. Let the eigenvalues of $L(G_i)$ be $\lambda_{i,1}\geq\cdots\geq\lambda_{i,n_i-1}>\lambda_{i,n_i}=0$. Then, for $i\in[k]$, there exist an orthogonal matrix $M_i$ such that $L(G_i)=M_iD_iM_i^*$ and $L(G)=MDM^*$ where $D_i=diag(\lambda_{i,1},\ldots,\lambda_{i,n_i-1},0)$, $M=M_1\oplus\cdots\oplus M_{k}$ and $D=D_1\oplus\cdots\oplus D_{k}$. Since any two $L_G(n,n-k)$-frames are unitary equivalent, there exist a unitary matrix $U$ such that  the synthesis operator of the $L_G(n,n-k)$-frame $F$ is $UB$, where $B=B_1\oplus B_2\oplus\cdots\oplus B_k$ such that $B_i=diag(\sqrt{\lambda_{i,1}},\ldots,\sqrt{\lambda_{i,n_i-1}})\widetilde{M_i}^*$ is the synthesis operator of the $L_{G_i}(n_i,n_i-1)$-frame $F_i=\{B_i(e_{i,j})\}_{j\in[n_i]}$ where $\widetilde{M_i}$ is formed by taking first $n_i-1$ columns of $M_i$ and $\{e_{i,j}\}_{j\in[n_i]}$ is the standard canonical orthonormal basis of $\mathbb{C}^{n_i}$. 
 
 Since $G_i$ is connected, by Theorem \ref{thm4.1}, $F_i$ is a full spark frame. Hence, for any $i\in[k]$, any $m-1$ columns of $B_i$ and thus of $B$ are linearly independent. Hence, any $ m-1 $ columns of $ UB $ are linearly independent. Let $m=n_r$ for some $r\in[k]$. Now $F_r$ is linearly dependent as $F_r$ is a frame for $\mathbb{C}^{n_r-1}$. Hence, the $m$ columns of $B_r$ are linearly dependent. Then, the corresponding $m$ columns of $B$  and thus of $ UB $ are linearly dependent. Therefore, the spark of the frame $F$ is $m$. 	\endproof

	\begin{thm}\label{thm4.4}
		Let $G$ be a graph of $n$ vertices with $k$ components and $F=\{f_i\}_{i\in[n]}$ be a $G(n,n-k)$-frame for $\mathbb{C}^{n-k}$. If $ m $ is the minimum cardinality of the vertex set of components of $ G $, then the spark of $F$ is $m$.
	\end{thm}
	\proof
	The proof follows from Theorem \ref{thm4.3} and Theorem \ref{thm2.6}.
	\endproof

	\section{Optimal dual frames of frames generated by graphs}\label{sec2}
	The idea of optimal dual frames (or simply $OD$-frame) was introduced by Lopez and Han \cite{LH1} to overcome the problem of erasures that arise from either transmission losses or disturbances during transmission. This section studies the $OD$-frames of frames generated by graphs. Let $F=\{f_i\}_{i\in[n]}$ be a frame. If $F$ has the frame operator $S$ and $c=\max\{\|f_i\|\|S^{-1}f_i\|:i\in[n]\}$, then $\Lambda_{1, F}:=\{l:\|f_l\|\|S^{-1}f_l\|=c\}$. We begin by showing that if the vertex set of a walk-regular component of the graph $G$ intersects $\Lambda_{1, F}$, then the canonical dual of an $L_G(n,k)$-frame $ F $ must be an $OD$-frame for $1$-erasure. 

\begin{thm}\label{thm5.1}
	Let $G$ be a graph with vertex set $ \{ 1,2, \dots, n\} $, $ k $ components, and has atleast a walk-regular component $ G_j $. If $F=\{f_i\}_{i\in[n]}$ is an $L_G(n,n-k)$-frame and $\Lambda_{1,F}\cap V(G_j)\neq\emptyset$, then the canonical dual frame of $F$ is an $OD$-frame for $1$-erasure.
\end{thm}

	\proof
	
	Let $G_1,\ldots,G_k$ be the components of $G$. Let $ |V(G_i)|=n_i $, for $ i\in [k] $ and $ l_j=n_1+\cdots+n_{j} $ for $ j\in [k]$. Suppose the vertex set of $G_i$ is $V(G_i)=\{l_{i-1}+1,l_{i-1}+2,\ldots,l_i\}$ where $l_0=0$. For $i\in [k]$, let $\lambda_{i,1}\geq\lambda_{i,2}\geq\cdots\geq\lambda_{i,n_i-1}>\lambda_{i,n_i}=0$  be the eigenvalues of $L(G_i)$. , there exist orthogonal matrices $M_i$ such that $L(G_i)=M_iD_iM_i^*$ where $D_i=diag(\lambda_{i,1},\lambda_{i,2},\ldots,\lambda_{i, n_i-1},0)$, for $i\in[k]$, and $L(G)=MDM^*$ where $M=M_1\oplus\cdots\oplus M_k$ and $D=D_1\oplus\cdots\oplus D_k$. Then, the synthesis operator of an $L_G(n,n-k)$-frame, say $\widetilde{F}=\{ \widetilde{f}_i\}$ is $B=B_1\oplus B_2\oplus\cdots\oplus B_k$ such that  $B_i=diag(\sqrt{\lambda_{i,1}},\sqrt{\lambda_{i,2}},\ldots,\sqrt{\lambda_{i,n_i-1}})\widetilde{M_i}^*$ (where $\widetilde{M_i}^*$ contains the first $n_i-1$ columns of $M_i$) is the synthesis operator of the $L_{G_i}(n_i,n_i-1)$-frame, say $F_i=\{f_{i,p}\}_{p\in[n_i]}$, for $\mathbb{C}^{n_i-1}$, for $i\in[k]$. By Theorem \ref{thm2.6}, there is a unitary operator $ U $ such that $ f_i=U\widetilde{f}_i $.

	Without loss of generality, assume $j=1$ that is $G_1$ is walk-regular and $\Lambda_{1,F}\cap V(G_1)\neq\emptyset$. Let $S$ and $\widetilde{S}$  be the frame operators of $F$  and $\widetilde{F}$, respectively. For $i\in [k]$, let $S_i$ be the frame operator of $F_i$. Then
	\begin{align*}
		\widetilde{S}&=BB^*\\
		&=(B_1\oplus\cdots\oplus B_k)(B_1^*\oplus\cdots\oplus B_k^*)\\
		&=B_1B_1^*\oplus\cdots\oplus B_kB_k^*\\
		&=S_1\oplus\cdots\oplus S_k	
	\end{align*}
	Thus, $\widetilde{S}^{-1}=S_1^{-1}\oplus\cdots\oplus S_k^{-1}$, and $ S^{-1}=U\widetilde{S}^{-1}U^* $. Hence, $S^{-1}(\widetilde{f}_i)=\left[\begin{array}{c}
		S_1^{-1}f_{1,i}\\
		\boldsymbol{0}
	\end{array}\right]$ for $i\in [n_1]$. 
	
	Since $G_1$ is a walk-regular graph, by Theorem \ref{thm3.3}, the Moore-Penrose inverse $L_1^+=M_1diag\left(\frac{1}{\lambda_{1,1}},\ldots,\frac{1}{\lambda_{1,n_1-1}},0\right)M_1^*$ has equal diagonal entries. Thus, there exists an $\alpha$ such that all the diagonal entries of  $L_1^+$ are equal to $\alpha$. By Theorem \ref{thm2.7}, $S_1^{-1}=diag\left(\frac{1}{\lambda_{1,1}},\ldots,\frac{1}{\lambda_{1,n_1-1}}\right)$. If $\{e_i\}_{i\in[n_1]}$ is the standard canonical orthonormal basis of $\mathbb{C}^{n_1}$, then $S_1^{-1}(f_{1,i})=S_1^{-1}B_1(e_i)=diag\left(\frac{1}{\sqrt{\lambda_{1,1}}},\ldots,\frac{1}{\sqrt{\lambda_{1,n_1-1}}}\right)\widetilde{M_1}^*(e_i)$ for $i\in[n_1]$. Thus, the synthesis operator of $\{S_1^{-1}f_{1,i}\}_{i\in[n_1]}$ is $T^*=diag\left(\frac{1}{\sqrt{\lambda_{1,1}}},\ldots,\frac{1}{\sqrt{\lambda_{1,n_1-1}}}\right)\widetilde{M_1}^*$. Hence, the Gramian matrix of $\{S_1^{-1}f_{1,i}\}_{i\in[n_1]}$ is $$TT^*=\widetilde{M_1}diag\left( \frac{1}{\lambda_{1,1}},\ldots,\frac{1}{\lambda_{1,n_1-1}}\right)\widetilde{M_1}^*=M_1diag\left(\frac{1}{\lambda_{1,1}},\ldots,\frac{1}{\lambda_{1,n_1-1}},0\right)M_1^*=L_1^+$$. Thus, for any $i\in[n_1]$, we have
	\begin{align*}
	\|S^{-1}(f_i)\|^2=\|U\widetilde{S}^{-1}U^*f_i\|^2=	\|\widetilde{S}^{-1}(\widetilde{f}_i)\|^2=\|S_1^{-1}f_{1,i}\|^2=(i,i)\text{-th entry of }L_1^+=\alpha.
	\end{align*}
	Since $G_1$ is a walk-regular graph, $G_1$ is a regular graph, say $r$-regular. Then, for any $i\in[n_1]$, we have
	$$\|S_1^{-1}f_{1,i}\|\|f_{1,i}\|=\sqrt{\alpha r}.$$
	Thus, by Theorem \ref{thm2.1}, the canonical dual $\{S_1^{-1}f_{1,i}\}_{i\in[n_1]}$ of the frame $\{f_{1,i}\}_{i\in[n_1]}$ is the unique $OD$-frame for any erasure. 
	
	Since $\|S^{-1}f_i\|\|f_i\|=\|\widetilde{S}^{-1}\widetilde{f}_i \| \|U \widetilde{f}_i \| = \|\widetilde{S}^{-1}\widetilde{f}_i \| \| \widetilde{f}_i \|=\|S_1^{-1}f_{1,i}\|\|f_{1,i}\|=\sqrt{\alpha r}$  for all $i\in[n_1]$ and $\Lambda_{1,F}\cap[n_1]\neq\emptyset$, $\max\{\|S^{-1}f_i\|\|f_i\|:i\in[n]\}=\sqrt{\alpha r}$.
	
	Suppose $H=\{h_i\}_{i\in[n]}$ be any dual frame of $F$. Then, $ U^*H=\{ U^*h_i\}_{i\in [n]} $ is a dual frame of $ \widetilde{F} $. Let $\widetilde{h}_i$ be obtained from $Uh_i$ by considering the first $n_1-1$ components, for $i\in[n_1]$. For any $f\in\mathbb{C}^{n_1-1}$ and $\widetilde{f}=\left[\begin{array}{c}
		f\\
		\boldsymbol{0}
	\end{array}\right]\in\mathbb{C}^{n-k}$, we have
	\begin{align*}
		\left[\begin{array}{c}
			f\\
			\boldsymbol{0}
		\end{array}\right]&=\sum_{i\in [n_1]}\langle \widetilde{f},U^*h_i\rangle \widetilde{f}_i\\
		&=\sum_{i\in[n_1]}\langle \widetilde{f},U^{*}h_i\rangle \left[\begin{array}{c}
			f_{1,i}\\
			\boldsymbol{0}
		\end{array}\right]\\
		&=\sum_{i\in [n_1]}\langle f,\widetilde{h}_i\rangle \left[\begin{array}{c}
			f_{1,i}\\
			\boldsymbol{0}
		\end{array}\right].
	\end{align*}
	Thus, $f=\sum_{i\in [n_1]}\langle f,\widetilde{h}_i\rangle f_{1,i}$ for any $f\in\mathbb{C}^{n_1-1}$. Thus, $\{\widetilde{h}_i\}_{i\in[n_1]}$ is a dual frame of $\{f_{1,i}\}_{i\in[n_1]}$. Since the canonical dual frame $\{S_1^{-1}f_{1,i}\}_{i\in[n_1]}$ of the frame $\{f_{1,i}\}_{i\in[n_1]}$ is the $OD$-frame of $\{f_{1,i}\}_{i\in[n_1]}$ for $1$-erasure, thus, 
	$$\max\{\|\widetilde{h}_i\|\|f_{1,i}\|:i\in[n_1]\}\geq\max\{\|S_1^{-1}f_{1,i}\|\|f_{1,i}\|:i\in[n]\}=\sqrt{\alpha r}.$$ 
  For any $i\in[n_1]$, $\|\widetilde{h}_i\|\leq\|U^*h_i\|=\|h_i\|$ and $\|f_i\|= \|U\widetilde{f}_i\|=\|\widetilde{f}_i \|=\|f_{1,i}\|$, thus, we have
  
	\begin{align*}
		\max\{\|h_i\|\|f_i\|:i\in[n]\}&\geq\max\{\|h_i\|\|f_i\|:i\in[n_1]\}\\
		&\geq\max\{\|\widetilde{h}_i\|\|f_{1,i}\|:i\in[n_1]\}\\
		&\geq \sqrt{\alpha r}\\
		&=\max\{\|S^{-1}f_i\|\|f_i\|:i\in[n]\}.
	\end{align*}
Hence, $D^1_{F,H}\geq D^1_{F,S^{-1}F}$ for any dual frame $H$ of $F$. Therefore, the canonical dual frame $S^{-1}F$ of $F$ is an $OD$-frame for $1$-erasure. 
	\endproof
	In the following corollary, we show that the Theorem \ref{thm5.1} also holds for the frames generated by graphs.

\begin{corollary}\label{cor5.1}
Let $G$ be a graph with vertex set $ \{ 1,2, \dots, n\} $ and has atleast a walk-regular component $ G_j $. If $F=\{f_i\}_{i\in[n]}$ is a frame generated by $ G $ and $\Lambda_{1, F}\cap V(G_j)\neq\emptyset$, then the canonical dual of $F$ is an $OD$-frame for $1$-erasure.
\end{corollary}
	\proof
	Suppose $G_1,\ldots,G_k$ are the components of $G$. Let $\widetilde{F}=\{\widetilde{f}_i\}_{i\in[n]}$ be an $L_G(n,n-k)$-frame for $\mathbb{C}^{n-k}$ 
 and $\widetilde{S}$ be the frame operator of $\widetilde{F}$. By Theorem \ref{thm2.6}, there is a unitary operator $U:\mathbb{C}^{n-k}\rightarrow\mathbb{C}^{n-k}$ such that $f_i=U(\widetilde{f}_i)$ for $i\in[n]$. If $S$ is the frame operator of $F$, then $S=U\widetilde{S}U^*$. Let $H=\{h_i\}_{i\in[n]}$ be any dual frame of $F$. Then, $U^*H=\{U^*h_i\}_{i\in[n]}$ is a dual frame of $\widetilde{F}$. Since the component $G_j$ is walk-regular and $\Lambda_{1,F}\cap V(G_j)\neq\emptyset$, by Theorem \ref{thm5.1}, $\widetilde{S}^{-1}\widetilde{F}$ is an $OD$-frame of $\widetilde{F}$ for $1$-erasure. Then, we have
	\begin{align*}
	\max\{\|h_i\|\|f_i\|:i\in[n]\}&=\max\{\|h_i\|\|U\widetilde{f}_i\|:i\in[n]\}\\
		&=\max\{\|U^*h_i\|\|\widetilde{f}_i\|:i\in[n]\}\\
		&\geq\max\{\|{\widetilde{S}}^{-1}\widetilde{f}_i\|\|\widetilde{f}_i\|:i\in[n]\}\\
		&=\max\{\|U^*S^{-1}U\widetilde{f}_i\|\|U\widetilde{f}_i\|:i\in[n]\}\\
		&=\max\{\|S^{-1}f_i\|\|f_i\|:i\in[n]\}.
	\end{align*}
	Thus, the canonical dual frame $S^{-1}F$ of $F$ is an $OD$-frame for $1$-erasure.
	\endproof
	
After studying the above results, the natural question arises whether the canonical dual is the unique $OD$-frame of the frame satisfying the conditions of Theorem \ref{thm5.1} and Corollary \ref{cor5.1}. In the next example, we present a frame generated by a graph satisfying the conditions of Theorem \ref{thm5.1}, but the canonical dual is not the unique $OD$-frame.
	\begin{example}{\em
		Consider the graph given in Figure \ref{fig1} of Example \ref{Exa4.1}. Then \begin{align*}
			F=\{f_i\}_{i\in [7]}=\left\{\left[\begin{array}{c}
				0\\
				0\\
				0\\
				\frac{\sqrt{6}}{2}\\
				\frac{-\sqrt{18}}{6}
			\end{array}\right],\left[\begin{array}{c}
				0\\
				0\\
				0\\
				0\\
				\frac{\sqrt{18}
				}{3}
			\end{array}\right], \left[\begin{array}{c}
				0\\
				0\\
				0\\
				\frac{-\sqrt{6}}{2}\\
				\frac{-\sqrt{18}}{6}
			\end{array}\right], \left[\begin{array}{c}
				1\\
				1\\
				0\\
				0\\
				0
			\end{array}\right], \left[\begin{array}{c}
				-1\\
				0\\
				1\\
				0\\
				0
			\end{array}\right], \left[\begin{array}{c}
				1\\
				-1\\
				0\\
				0\\
				0
			\end{array}\right], \left[\begin{array}{c}
				-1\\
				0\\
				-1\\
				0\\
				0
			\end{array}\right]\right\}
		\end{align*} is an $L_G(7,5)$-frame for $\mathbb{C}^5$. Note that each connected component of $G$ is walk-regular.  Hence, by Theorem \ref{thm5.1}, $S^{-1}F$ is an $OD$-frame for $1$-erasure, where $S$ is the frame operator of $F$. Let $M_1$ be the same as in Example \ref{Exa4.1}. Then, the synthesis operator of $S^{-1}F$ is $diag(\frac{1}{2},\frac{1}{\sqrt{2}},\frac{1}{\sqrt{2}},\frac{1}{\sqrt{3}},\frac{1}{\sqrt{3}})M_1^*$. Thus, we have
		\begin{align*}
			S^{-1}F=\left\{\left[\begin{array}{c}
				0\\
				0\\
				0\\
				\frac{1}{\sqrt{6}}\\
				\frac{-1}{\sqrt{18}}
			\end{array}\right],\left[\begin{array}{c}
				0\\
				0\\
				0\\
				0\\
				\frac{\sqrt{18}
				}{9}
			\end{array}\right], \left[\begin{array}{c}
				0\\
				0\\
				0\\
				\frac{-1}{\sqrt{6}}\\
				\frac{-1}{\sqrt{18}}
			\end{array}\right], \left[\begin{array}{c}
				\frac{1}{4}\\
				\frac{1}{2}\\
				0\\
				0\\
				0
			\end{array}\right], \left[\begin{array}{c}
				\frac{-1}{4}\\
				0\\
				\frac{1}{2}\\
				0\\
				0
			\end{array}\right], \left[\begin{array}{c}
				\frac{1}{4}\\
				\frac{-1}{2}\\
				0\\
				0\\
				0
			\end{array}\right], \left[\begin{array}{c}
				\frac{-1}{4}\\
				0\\
				\frac{-1}{2}\\
				0\\
				0
			\end{array}\right]\right\}
		\end{align*}
		Now $\max\{\|S^{-1}f_i\|\|f_i\|:i\in[7]\}=\max\left\{\frac{2}{3},\frac{\sqrt{10}}{4}\right\}=\frac{\sqrt{10}}{4}$. Let $\mu=\left[\begin{array}{c}
			0\\
			0\\
			0\\
			0.01\\
			0.01
		\end{array}\right]\in\mathbb{C}^5$. Then, by Theorem \ref{thm2.8},  $H=\{h_i\}_{i\in[7]}=\{S^{-1}f_1+\mu,S^{-1}f_2+\mu,S^{-1}f_3+\mu,S^{-1}f_4,S^{-1}f_5,S^{-1}f_6,S^{-1}f_7\}$ is an alternate dual frame of $F$. Then, $\max\{\|h_i\|\|f_i\|:i\in[7]\}=\max\left\{0.672,0.647,0.681,\frac{\sqrt{10}}{4}\right\}\break= \frac{\sqrt{10}}{4}=\max\{\|S^{-1}f_i\|\|f_i\|:i\in[7]\}$. Hence, $H$ is also an $OD$-frame of $F$ for $1$-erasure. Thus, $S^{-1}F$ is not the unique $OD$-frame.}
	\end{example}
	
	Next, we prove that the canonical dual of a frame generated by a walk-regular graph is the unique $OD$-frame.
	\begin{thm}\label{thm5.2}
		If $G$ is a walk-regular graph and $F=\{f_i\}_{i\in[n]}$ is an $L_G(n,k)$-frame, then the canonical dual frame of $F$ is the unique $OD$-frame for any erasure.
	\end{thm}
	\proof
	Suppose $S$ is the frame operator of $F$. By Theorem \ref{thm3.3}, diagonal entries of the Moore-Penrose inverse $L^+$ of the Laplacian matrix $L$ are equal. Suppose the diagonal entries of $ L^+ $ is $ \alpha $. Then, using the same argument as in the proof of Theorem \ref{thm5.1}, $\|S^{-1}f_i\|^2= (i,i)\text{-th diagonal entry of }  L^+ =\alpha$. Since walk-regular graphs are regular, assume $G$ is an $r$-regular graph. Then, for $i\in[n]$,  $\|S^{-1}f_i\|\|f_i\|=\sqrt{\alpha r}$. Hence, by Theorem \ref{thm2.1}, the canonical dual of $F$ is the unique $OD$-frame for any erasure.
	\endproof
	
As a consequence of the above theorem, we obtain the following corollary.
	\begin{corollary}
		If $G$ is a walk-regular graph and $F=\{f_i\}_{i\in[n]}$ is a frame generated by $G$, then the canonical dual frame of $F$ is the unique $OD$-frame for any erasure.
	\end{corollary}

One may think whether the results proved in this section hold for regular graphs, which are not walk-regular. The answer is negative. In Example \ref{ex5.2}, we show a frame generated by a connected regular graph that is not walk-regular such that its canonical dual is not an $OD$-frame for $1$-erasure.
	\begin{example} \label{ex5.2}{\em
		Consider the graph $ G $ given in Figure \ref{fig2}. One may note that $G$ is a regular graph but not a walk-regular graph. Now, we generate an $L_G(8,7)$-frame for $\mathbb{C}^7$.  Then  
		
		\begin{figure}
			\begin{center}
				\includegraphics[scale= 0.65]{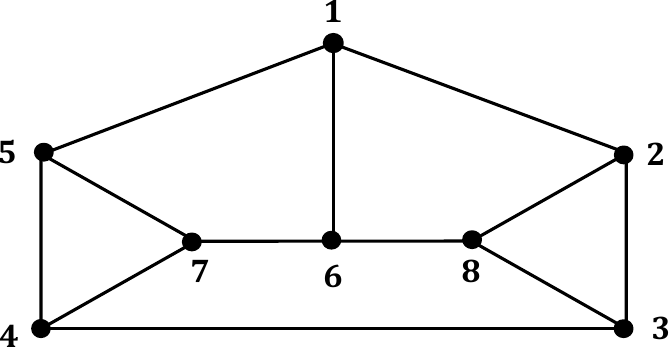}
				\caption{Graph $ G $} \label{fig2}
			\end{center}
		\end{figure}

		\begin{align*}
			L(G)=	\left[\begin{array}{cccccccc}
				3 & -1 & 0 & 0 & -1 & -1 & 0 & 0\\
				-1 & 3 & -1 & 0 & 0 & 0 & 0 & -1\\
				0 & -1 & 3 & -1 & 0 & 0 & 0 & -1\\
				0 & 0 & -1 & 3 & -1 & 0 & -1 & 0\\
				-1 & 0 & 0 & -1 & 3  & 0 & -1 & 0\\
				-1 & 0 & 0 & 0 & 0 & 3 & -1 & -1\\
				0 & 0 & 0 & -1 & -1 & -1 & 3 & 0\\
				0 & -1 & -1 & 0 & 0 & -1 & 0 & 3
			\end{array} 
			\right]=MDM^*
		\end{align*}
		where 
		\begin{align*}
			M=	\left[\begin{array}{cccccccc}
				\frac{\sqrt{3}}{6} & \frac{-\sqrt{6}}{12} & \frac{1}{2} & \frac{1}{2} & \frac{1}{2} & 0 & 0 & \frac{\sqrt{2}}{4} \\
				0 & \frac{\sqrt{6}}{4} & 0 & \frac{\sqrt{2}}{4} & \frac{-\sqrt{2}}{4} & \frac{1}{2\sqrt{3-\sqrt{3}}} & \frac{1}{2\sqrt{3+\sqrt{3}}} & \frac{\sqrt{2}}{4} \\
				\frac{\sqrt{3}}{6} & \frac{-\sqrt{6}}{12} & \frac{-1}{2} & 0 & 0 & \frac{\sqrt{3}-1}{2\sqrt{3-\sqrt{3}}} & \frac{-(\sqrt{3}+1)}{2\sqrt{3+\sqrt{3}}} & \frac{\sqrt{2}}{4} \\
				\frac{\sqrt{3}}{6} & \frac{-\sqrt{6}}{12} & \frac{-1}{2} & 0 & 0 & \frac{1-\sqrt{3}}{2\sqrt{3-\sqrt{3}}} & \frac{\sqrt{3}+1}{2\sqrt{3+\sqrt{3}}} & \frac{\sqrt{2}}{4} \\
				-\frac{\sqrt{3}}{3} & \frac{-\sqrt{6}}{12} & 0 & \frac{\sqrt{2}}{4} & \frac{-\sqrt{2}}{4} & \frac{-1}{2\sqrt{3-\sqrt{3}}} & \frac{-1}{2\sqrt{3+\sqrt{3}}} & \frac{\sqrt{2}}{4} \\
				\frac{\sqrt{3}}{6} & \frac{-\sqrt{6}}{12} & \frac{1}{2} & \frac{-1}{2} & \frac{-1}{2} & 0 & 0 & \frac{\sqrt{2}}{4} \\
				0 & \frac{\sqrt{6}}{4} & 0 & \frac{-\sqrt{2}}{4} & \frac{\sqrt{2}}{4} & \frac{-1}{2\sqrt{3-\sqrt{3}}} & \frac{-1}{2\sqrt{3+\sqrt{3}}} & \frac{\sqrt{2}}{4} \\
				\frac{-\sqrt{3}}{3} & \frac{-\sqrt{6}}{12} & 0 & \frac{-\sqrt{2}}{4} & \frac{\sqrt{2}}{4} & \frac{1}{2\sqrt{3-\sqrt{3}}} & \frac{1}{2\sqrt{3+\sqrt{3}}} & \frac{\sqrt{2}}{4} 
			\end{array} 
			\right]
		\end{align*}
\newline
		and $D=diag(4,4,2,4-\sqrt{2},4+\sqrt{2},3-\sqrt{3},3+\sqrt{3},0)$. If $M_1$ is formed by taking the first seven columns of $M$, $D_1=diag(2,2,\sqrt{2},\sqrt{4-\sqrt{2}},\sqrt{4+\sqrt{2}},\sqrt{3-\sqrt{3}},\sqrt{3+\sqrt{3}})$ and $B=D_1M_1^*$, then $F=\{B(e_i)\}_{i\in[8]}$ is an $L_G(8,7)$-frame for $\mathbb{C}^7$. Also, $\|f_i\|^2=i$-th diagonal entry of $L(G)=3$ for all $i\in[8]$. Let $D_2=diag(\frac{1}{2},\frac{1}{2},\frac{1}{\sqrt{2}},\frac{1}{\sqrt{4-\sqrt{2}}},\frac{1}{\sqrt{4+\sqrt{2}}},\frac{1}{\sqrt{3-\sqrt{3}}},\frac{1}{\sqrt{3+\sqrt{3}}})$. Thus, the canonical dual frame of $F$ is
  \vspace{10pt}
	\begin{align*}
			S^{-1}F=D_2M_1^*= \left \{ \left[\begin{array}{c}
				\frac{\sqrt{3}}{12}\\
				\frac{-\sqrt{6}}{24}\\
				\frac{\sqrt{2}}{4}\\
				\frac{1}{2\sqrt{4-\sqrt{2}}}\\
				\frac{1}{2\sqrt{4+\sqrt{2}}}\\
				0\\
				0
			\end{array}\right],\left[\begin{array}{c}
				0\\
				\frac{\sqrt{6}}{8}\\
				0\\
				\frac{\sqrt{2}}{4\sqrt{4-\sqrt{2}}}\\
				\frac{-\sqrt{2}}{4\sqrt{4+\sqrt{2}}}\\
				\frac{-1}{2(\sqrt{3}-3)}\\
				\frac{1}{2(\sqrt{3}+3)}
			\end{array}\right], \left[\begin{array}{c}
				\frac{\sqrt{3}}{12}\\
				\frac{-\sqrt{6}}{24}\\
				\frac{-\sqrt{2}}{4}\\
				0\\
				0\\
				\frac{-(\sqrt{3}-1)}{2(\sqrt{3}-3)}\\
				\frac{-(\sqrt{3}+1)}{2(\sqrt{3}+3)}
			\end{array}\right], \left[\begin{array}{c}
				\frac{\sqrt{3}}{12}\\
				\frac{-\sqrt{6}}{24}\\
				\frac{-\sqrt{2}}{4}\\
				0\\
				0\\
				\frac{\sqrt{3}-1}{2(\sqrt{3}-3)}\\
				\frac{\sqrt{3}+1}{2(\sqrt{3}+3)}
			\end{array}\right], \left[\begin{array}{c}
				\frac{-\sqrt{3}}{6}\\
				\frac{-\sqrt{6}}{24}\\
				0\\
				\frac{\sqrt{2}}{4\sqrt{4-\sqrt{2}}}\\
				\frac{-\sqrt{2}}{4\sqrt{4+\sqrt{2}}}\\
				\frac{1}{2(\sqrt{3}-3)}\\
				\frac{-1}{2(\sqrt{3}+3)}
			\end{array}\right]\right., 
   \end{align*}

\begin{align*}
\hspace{2cm}\left.\left[\begin{array}{c}
				\frac{\sqrt{3}}{12}\\
				\frac{-\sqrt{6}}{24}\\
				\frac{\sqrt{2}}{4}\\
				\frac{-1}{2\sqrt{4-\sqrt{2}}}\\
				\frac{-1}{2\sqrt{4+\sqrt{2}}}\\
				0\\
				0
			\end{array}\right],
			\left[\begin{array}{c}
				0\\
				\frac{\sqrt{6}}{8}\\
				0\\
				\frac{-\sqrt{2}}{4\sqrt{4-\sqrt{2}}}\\
				\frac{\sqrt{2}}{4\sqrt{4+\sqrt{2}}}\\
				\frac{1}{2(\sqrt{3}-3)}\\
				\frac{-1}{2(\sqrt{3}+3)}
			\end{array}\right],\left[\begin{array}{c}
				\frac{-\sqrt{3}}{6}\\
				\frac{-\sqrt{6}}{24}\\
				0\\
				\frac{-\sqrt{2}}{4\sqrt{4-\sqrt{2}}}\\
				\frac{\sqrt{2}}{4\sqrt{4+\sqrt{2}}}\\
				\frac{-1}{2(\sqrt{3}-3)}\\
				\frac{1}{2(\sqrt{3}+3)}
			\end{array}\right]\right \}
		\end{align*}


		Then, $\|S^{-1}f_1\|=\|S^{-1}f_6\|=0.5469$, $\|S^{-1}f_2\|=\|S^{-1}f_5\|=\|S^{-1}f_7\|=\|S^{-1}f_8\|=0.5761$, and $\|S^{-1}f_3\|=\|S^{-1}f_4\|=0.5682$. Thus, we have
		\begin{align*}
		D^1_{F,S^{-1}F}=\max\{\|S^{-1}f_i\|\|f_i\|:i\in[8]\}=0.9978.
		\end{align*}
		Take $\mu=[0.001,-0.001,0,0,0,0,0,0]^t\in\mathbb{C}^8$. Then, by Theorem \ref{thm2.8}, $H=\{h_i\}_{i\in[8]}=\{ S^{-1}f_i+\mu\}_{i\in [8]}$ is a dual frame of $F$ such that 
		\begin{align*}
			D^1_{F,H}=\max\{\|h_i\|\|f_i\|:i\in[8]\}&=\sqrt{3}\max\{0.5474,0.5755,0.5687,0.5757\}\\
			&=0.9971\\
			&<D^1_{F,S^{-1}F}.
		\end{align*}
Thus,  $S^{-1}F$ is not an $OD$-frame for $1$-erasure.}
	\end{example}
In Theorem \ref{thm2.1}, Lopez and Han provide a sufficient condition under which the canonical dual is the unique $OD$-frame. In the next theorem, we establish that the condition is also necessary for frames generated by connected graphs.
	\begin{thm}\label{thm1}
		Let $G$ be a connected graph and $F=\{f_i\}_{i\in[n]}$ be a $G(n,n-1)$-frame for $\mathbb{C}^{n-1}$ with frame operator $S$. Then, $\{S^{-1}f_i\}_{i\in[n]}$ is the unique $OD$-frame for any erasure if and only if  $\|S^{-1}f_i\|\|f_i\|$ is constant for all $i\in[n]$.
	\end{thm}
	\proof
	First suppose that $\{S^{-1}f_i\}_{i\in[n]}$ is the unique $OD$-frame of $F$ for any erasure. Assume that $\|S^{-1}f_i\|\|f_i\|$ is not constant for $i\in[n]$. Then, $|\Lambda_{1,F}|\leq n-1$.
		
	It is given that $G$ is a connected graph, thus, by Theorem \ref{thm4.2}, $F$ is a full spark frame for $\mathbb{C}^{n-1}$. Hence, $\{f_i\}_{i\in\Lambda_{1,F}}$ is linearly independent. 
	 Also, there exist non-zero scalars $\alpha_1,\alpha_2,\ldots,\alpha_n$ such that $\sum\limits_{i=1}^n\alpha_if_i=0$. Thus, by Theorem \ref{thm2.2}, $\{S^{-1}f_i\}_{i\in[n]}$ is not an $OD$-frame for $1$-erasure, which is a contradiction. Thus, $\|S^{-1}f_i\|\|f_i\|$ is constant for all $i\in[n]$.
	 
	Converse follows straight forward from Theorem \ref{thm2.1}.
	\endproof
	We end this article by showing that the canonical dual frames of frames generated by connected graphs are either the unique $OD$-frame or not an $OD$-frame for $1$-erasure.
	\begin{thm}
		Let $G$ be a connected graph and $F=\{f_i\}_{i\in[n]}$ be a frame generated by the graph $G$ with frame operator $S$. Then, either $\{S^{-1}f_i\}_{i\in[n]}$ is the unique $OD$-frame for $1$-erasure or $\{S^{-1}f_i\}_{i\in[n]}$ is not an $OD$-frame for $1$-erasure.
	\end{thm}
	\proof
	If $\{S^{-1}f_i\}_{i\in[n]}$ is the unique $OD$-frame for $1$-erasure, then we are done. Now assume that $\{S^{-1}f_i\}_{i\in[n]}$ is not the unique $OD$-frame for $1$-erasure. Then, by Theorem \ref{thm1}, $|\Lambda_{1,F}|\leq n-1$. Also, by Theorem \ref{thm4.2}, $F$ is a full spark frame for $\mathbb{C}^{n-1}$. Thus, $\{f_i\}_{i\in\Lambda_{1,F}}$ is linearly independent and there exist non-zero scalars $\alpha_1,\alpha_2,\ldots,\alpha_n$ such that $\sum\limits_{i\in[n]}\alpha_if_i=0$. Thus, by Theorem \ref{thm2.2}, $\{S^{-1}f_i\}_{i\in[n]}$ not an $OD$-frame for $1$-erasure.
	\endproof
	
	\section*{Acknowledgments}
	Aniruddha Samanta expresses thanks to the National Board for Higher Mathematics (NBHM), Department of Atomic Energy, India, for providing financial support in the form of an NBHM Post-doctoral Fellowship (Sanction Order No. 0204/21/2023/R\&D-II/10038). The second author also acknowledges excellent working conditions in the Theoretical Statistics and Mathematics Unit, Indian Statistical Institute Kolkata. 

	\mbox{}
		
\end{document}